\documentclass{amsart}

\usepackage{amstext,amssymb,amsmath,amsbsy,mathrsfs,mathtools,amsfonts,amscd,exscale}
\usepackage{graphics, epstopdf}
\usepackage{color}
\usepackage{amsmath, amssymb}
\usepackage{enumerate}
\usepackage{bigints}
\usepackage{verbatim}
\usepackage[poorman]{cleveref}
\usepackage{txfonts} 
\usepackage[normalem]{ulem}

\usepackage{graphicx}
\usepackage{epstopdf}
\usepackage{float}
\usepackage{placeins}
\usepackage{float}
\usepackage[caption = false]{subfig}
\usepackage{multirow}

\newcommand{\calE}{\mathcal{E}}
\newcommand{\calF}{\mathcal{F}}

\renewcommand{\sc}{\textsc}

\newcommand{\bnu}{{\boldsymbol \nu}}

\usepackage{algorithm}
\usepackage{algpseudocode}
\newcommand{\algrule}[1][.2pt]{\par\vskip.5\baselineskip\hrule height #1\par\vskip.5\baselineskip}

\theoremstyle{definition}

\theoremstyle{remark}

\numberwithin{equation}{section}

\newcommand{\abs}[1]{\lvert#1\rvert}

\newcommand\modif[1]{{#1}}

\def\l({\left(}
\def\r){\right)}

\def\nd{\partial_{\mu}}

\def\LtGb{L^2(\Gamma_h^b)}
\def\LtGi{L^2(\Gamma_h^0)}

\def\Vh{\mathbb{V}_h }
\def\Vhk{\mathbb{V}_h^k }

\def\lom{{L^2(\Omega)}}

\def\htwom{{H^2(\Omega)}}

\def\o{\Omega}
\newcommand{\Th}{\mathcal{T}_h}

\def\Ei{\mathcal{E}_h^0}
\def\Eb{\mathcal{E}_h^b}

\def\Ahk{\mathbb{A}_{h,\delta}} 

\definecolor{red}{rgb}{1,0,0}
\definecolor{blue}{rgb}{0,0,1}

\definecolor{aggiemaroon}{rgb}{0.8,0,0}

\begin{document}
	
\title[Large Bending Deformation of a Bilayer Plate with Isometry Constraint]{Discontinuous Galerkin Approach to Large Bending Deformation of a Bilayer Plate with Isometry Constraint}

      \author{Andrea Bonito$^1$}
        \address{Department of Mathematics, Texas A\&M University, College Station, TX 77843}
        \email{bonito@math.tamu.edu}
        \thanks{$^1$ Partially supported by the NSF Grant DMS-1817691}
        
        \author{Ricardo H. Nochetto$^2$}
        \address{Department of Mathematics and Institute for Physical Science
                and Technology, University of Maryland, College Park, Maryland 20742}
        \email{rhn@math.umd.edu}
        \thanks{$^2$ Partially supported by the NSF Grants DMS-1411808 and DMS-1908267, the
                Institut Henri Poincar\'e (Paris) and the Hausdorff Institute (Bonn).}

        \author{Dimitris Ntogkas$^3$}
        \address{Wells Fargo, North Carolina 28202}
        \email{dimitnt@gmail.com}
        \thanks{$^3$ Partially supported by the  NSF Grant DMS-1411808 and
                the 2016-2017 Patrick and Marguerite Sung Fellowship of the
                University of Maryland.}
\begin{abstract}
We present a computational model of thin elastic bilayers that undergo large bending isometric deformations when actuated by non-mechanical stimuli. We propose a discontinuous Galerkin approximation of the variational formulation discussed in \cite{BaBoMuNo}. We showcase its advantages and good computational performance with configurations of interest in both engineering and medicine and either Dirichlet or free boundary conditions.
\end{abstract}

\maketitle

\section{Introduction} \label{S:Introduction}

Large bending deformation of thin plates is a critical subject for many modern engineering and medical
applications due to the extensive use of plate actuators in a variety of systems like thermostats, nano-tubes, micro-robots and micro-capsules \cite{BasAbLaGr,JaSmIn,KuLPL,SchmidtEb,SmInLu}. From the mathematical viewpoint, there is an increasing interest in the modeling and the numerical treatment of such plates.
A rigorous analysis of large bending deformation of plates was conducted in the seminal work of Friesecke, James and M\"uller \cite{FJM}, who derive
geometrically non-linear Kirchhoff models from three dimensional hyperelasticity via $\Gamma$-convergence. There have been various other interesting models since then, such as the prestrained model derived in \cite{Lewicka}. 
A well-known relevant case that provides also the main motivation for this work is the bilayer-plate model \cite{BaBoMuNo,BaBoNo,Schmidt}. This is the case when two thin layers with different material properties are attached together. Upon thermal, electrical, or chemical actuation, such plates react differently thereby leading to large bending deformations. Mathematically, the plates develop an intrinsic spontaneous curvature tensor $Z$ and the deformation $y:\Omega\to\mathbb{R}^3$ of the midplane $\Omega\subset\mathbb{R}^2$ minimizes the elastic energy
\begin{equation}\label{E:bilayer_energy}
    E[y] := \frac{1}{2} \int_{\Omega} |H - Z|^2 - \int_{\Omega} f y,
\end{equation}  
where $H = (h_{ij})_{i,j=1}^2$ is the {\it second fundamental form} of the surface $y(\Omega)$, namely $h_{ij} = \partial_{ij} y \cdot \bnu$, $\bnu$ is the unit normal to the midplate and $|\cdot|$ denotes the Frobenius norm.
This spontaneous deformation of a flat configuration may occur without a forcing term $f$. Our
goal is to explore simple mathematical and computational bilayer models which are capable to reproduce equilibrium configurations observed in engineering and medical applications \cite{AlBSm,JHZ,Simpson}. Shear and stretch are negligible leaving \emph{bending} as the essential feature to account for.
In particular, lengths and areas are preserved by the admissible deformations $y(\Omega)$.
Mathematically, this entails that $y$ is an {\it isometry}, namely $y$ obeys the pointwise constraint
\begin{equation}\label{E:isometry}
[\nabla y]^{\mathsf{T}} \nabla y = I
\end{equation}
for the {\it first fundamental form},
where $\nabla y = [\partial_1 y, \partial_2 y]\in\mathbb{R}^{3\times2}$, $[\nabla y]^{\mathsf{T}}$ is the transpose of $\nabla y$, and $I$ is the $2\times2$ identity matrix. For isometric deformations, the unit normal $\bnu$ reads
\begin{equation}\label{E:normal}
\bnu := \partial_1 y \times \partial_2 y.
\end{equation}    
\modif{Before proceeding further, we point out that because of the isometry constraint \eqref{E:isometry} (no shear nor stretch),  the Lam\'e's first parameter $\lambda$ and shear-modulus $\mu$ of the plate only affect the energy by a multiplicative constant $\frac{1}{6} \frac{2\mu +\lambda \mu}{2\mu+\lambda}$ which is thus omitted \cite{FJM,Schmidt}.}

Previous work on the numerical treatment of large bending deformations includes the use of Kirchhoff finite element discretizations for the single layer problem by Bartels \cite{Bartels} and the bilayer problem by Bartels, Bonito and Nochetto \cite{BaBoNo}, both reviewed by Bartels \cite{Bartels-HB}. Modeling and simulation of thermally actuated bilayer plates is developed by Bartels, Bonito, Muliana and Nochetto \cite{BaBoMuNo}. Our discontinuous Galerkin (dG) approach for the single layer problem, described in \cite{BoNoNt}, exhibits some desirable theoretical and computational properties. We show $\Gamma$-convergence in \cite{BoNoNt} of the discrete energy functional $E_h$ to the continuous bending energy $E$ and combine this property with a compactness argument to deduce that global minimizers of $E_h$ converge strongly (up to a subsequence) to global minimizers of $E$ in $\lom$. We also test the dG method computationally in \cite{BoNoNt} via some illuminating examples that document its improved accuracy and geometric flexibility with respect to the Kirchhoff approach in \cite{Bartels}. These experiments are simple but designed to explore theoretical properties of the dG method and showcase its potential.

In this work we extend the dG approach to the simplified bilayer model of \cite{BaBoMuNo} and use it to compute configurations of interest in engineering and medical applications with either Dirichlet or free boundary conditions. In Section \ref{S:SingleLayer} we recall the single layer model of \cite{BoNoNt} and its dG approximation. In Section \ref{S:SimpleBilayer}, we briefly discuss the bilayer model of \cite{BaBoMuNo} and its novel dG approximation. We conclude in Section \ref{S:NumericalExperiments} with simulations relevant in applications. Some simulations reproduce those already presented in \cite{BaBoNo,BaBoMuNo} using Kirchhoff elements and some are new. 

\section{Single Layer - Main features of the DG Approach} \label{S:SingleLayer}
We start with a brief description of the single layer framework, namely $Z=0$ in \eqref{E:bilayer_energy}, and the key features of the dG method developed in \cite{BoNoNt}.

\medskip
{\bf Continuous energy.}
Given a Lipschitz domain $\Omega \subset \mathbb R^2$ (the undeformed plate), a smooth isometric deformation $y: \Omega \to \mathbb{R}^3$ satisfies
\begin{equation}\label{e:frob}
|H|^2 = |D^2 y|^2 = |\Delta y|^2
\end{equation}
pointwise in $\Omega$. To see this, we write the isometry constraint~\eqref{E:isometry} componentwise
 $$
 \partial_i y \cdot \partial_j y = \delta_{ij}, \quad i,j=1,2,
$$
and deduce by differentiation
$$
\partial_k y \cdot \partial_{ij} y = 0
\qquad\textrm{ or }\qquad
\partial_{ij} y \parallel \bnu = \partial_1 y \times \partial_2 y, \quad i,j,k=1,2.
$$
Combining this with the definition of the second fundamental form $H=(h_{ij})_{ij=1}^2$, leads to 
\begin{equation}\label{E:hessian}
h_{ij} = \partial_{ij} y \cdot \bnu
\qquad\Rightarrow\qquad
\partial_{ij} y =  h_{ij} \bnu = h_{ij} \, \partial_1 y \times \partial_2 y,
\end{equation}
which implies the first equality in \eqref{e:frob}. The second one follows from
\[
\partial_1(\partial_{12} y\cdot\partial_2 y)
= \partial_2(\partial_{11} y \cdot\partial_2 y)  = 0
\quad\Rightarrow\quad
|\partial_{12} y|^2 = \partial_{11} y \cdot \partial_{22} y.
\]
Therefore, the elastic energy \eqref{E:bilayer_energy} for the single layer plate
can be rewritten as
\begin{equation}\label{E:energy-single-layer}
E^0[y] := \frac12 \int_\Omega |D^2 y|^2 - \int_\Omega fy,
\end{equation}
with load $f \in [\lom]^3$. The single layer model thus consists of seeking a
deformation $y \in [H^2(\Omega)]^3$ that minimizes \eqref{E:energy-single-layer}
subject to the nonlinear pointwise isometry constraint \eqref{E:isometry} a.e. in $\Omega$
as well as {\it Dirichlet} boundary conditions on
$\partial_D \Omega \subset  \partial  \Omega$
\begin{equation}\label{E:bc}
y=g, \quad \nabla y = \Phi \qquad\textrm{on }\partial_D \Omega,
\end{equation}
where $g \in [H^2(\Omega)]^3, \Phi:= \nabla g \in [H^1(\Omega)]^{3\times 2}$ are given and satisfy $\Phi^{\mathsf{T}} \Phi = I$ on $\partial_D\Omega$. In addition to \cite{BoNoNt}, we allow $\partial_D\Omega=\emptyset$ and call it {\it free} boundary conditions.

The first variation  of \eqref{E:isometry} yields the {\it linearized isometry constraint}
\begin{equation}\label{E:linear-isometry}
  L[v;y] := [\nabla v]^{\mathsf{T}} \nabla y + [\nabla y]^{\mathsf{T}} \nabla v = 0
  \quad\text{a.e. } \Omega.
\end{equation}
This defines the tangent plane $\calF[y]$ at $y$ to the nonlinear constraint \eqref{E:isometry}, namely the set of functions $v\in[H^2(\Omega)]^3$ with vanishing Dirichlet
trace on $\partial_D\Omega$ provided $\partial_D\Omega\ne\emptyset$
and satisfying \eqref{E:linear-isometry}.
The first variation $\delta E^0[y](v)=0$ of \eqref{E:energy-single-layer} in the direction $v$ reads
\begin{equation}\label{E:EL}
  \delta E^0[y](v) = \int_\Omega D^2 y : D^2 v \, - \, f\cdot v = 0
  \qquad\forall \, v\in\calF[y],
\end{equation}
and is the Euler-Lagrange equation. We can find its solution $y\in [H^2(\Omega)]^3$  as the limit as $t\to\infty$ of the solution $y(t)\in  [H^2(\Omega)]^3$ with $\partial_t y(t) \in \calF[y(t)]$ of the $H^2$-gradient flow
\begin{equation}\label{E:gf-1}
\big( \partial_t y, v \big)_{H^2(\Omega)} +\modif{\delta E^0[y](v)} = 0
\qquad\forall \, v\in \calF[y].
\end{equation}
For $v,w \in [H^2(\Omega)]^3$, the $H^2(\Omega)$-scalar product in \eqref{E:gf-1}
is defined to be
\begin{equation}\label{e:H2}
\big(v,w)_{H^2(\Omega)}:= \int_\Omega D^2 v : D^2 w + \varepsilon \int_\Omega v \cdot w.
\end{equation}
Hereafter, we take $\varepsilon >0$ whenever $\partial_D \Omega = \emptyset$ and $\varepsilon = 0$ provided $\partial_D \Omega \not =  \emptyset$.
The evolution \eqref{E:gf-1} is supplemented by an initial condition $y(0) = y_0$ where $y_0 \in [H^2(\Omega)]^3$ satisfies both \eqref{E:isometry} and \eqref{E:bc}.
Notice that this, together with $\partial_t y\in\calF[y]$, yields
$$
[\nabla y(t)]^{\mathsf{T}} \nabla y(t) - I = \int_0^t \partial_t \left([\nabla y(s)]^{\mathsf{T}} \nabla y(s)\right)ds = \int_0^t L[\partial_s y(s); y(s)] ds = 0,
$$
whence $y(t)$ remains an isometry along the gradient flow \eqref{E:gf-1}.
Moreover, since $\partial_t y\in\calF[y]$, taking $v=\partial_t y$ in \eqref{E:gf-1} gives
\begin{equation}\label{E:energy-decrease}
  \| \,\partial_t y \,\|_{H^2(\Omega)}^2 + \frac{d}{dt} E^0[y] = 0
  \qquad\Rightarrow\qquad \frac{d}{dt} E^0[y] \le 0.
\end{equation}

\medskip
{\bf Discrete energy.}
We consider  a sequence of subdivisions $\{ \mathcal T_h \}_{h>0}$ of $\Omega$  made of triangles or quadrilaterals.
We assume that the sequence is shape regular, quasi-uniform and identify $h$ with the maximal element size.
From now on $c$ and $C$ are generic constants independent of $h$ but possibly depending on the shape-regularity and quasi-uniformity constants of the sequence $\{ \mathcal T_h \}_{h >0}$.

We denote by $\mathbb{P}_k$ (resp. $\mathbb Q_k$) the space of polynomial functions of degree at most $k\geq 0$ (resp. at most $k$ on each variable). Also, $\widehat T$ stands for the reference element, which is either the master triangle when the subdivision is made of triangles or the unit square in the case of quadrilaterals. The mapping between the reference element $\widehat T$ and any $T \in \mathcal T_h$ is denoted $F_T$. Notice that $F_T$ is affine for triangles $T$ and bi-linear for quadrilaterals $T$. 

With each subdivision $\mathcal T_h$ made of triangles, we associate the space of discontinuous piecewise polynomial functions
\begin{equation}\label{E:discrete-space}
\Vhk := \left\{ v_h \in \lom \  :  \ \ \modif{v_h \circ  F_T  \in \mathbb{P}_k, \quad \forall T \in \mathcal T_h} \right\}.
\end{equation}
Alternatively for subdivisions made of quadrilaterals, 
the space $\mathbb{P}_k$ is replaced by  $\mathbb{Q}_k$. We point out that in this case, $v_h|_T$ is no longer polynomial, which entails additional difficulties in the analysis. We refer to \cite{BoNoNt} for details but note that we require $k\geq 2$ in both cases.

We denote by $\Ei$ the collection of edges of $\mathcal T_h$ contained in $\Omega$ and by $\Eb$ those contained in $\partial_D\Omega$ (note that $\Eb=\emptyset$ provided $\partial_D\Omega=\emptyset$); hence $\calE_h:=\Ei\cup\Eb$ is the set of active interelement boundaries (across which jumps and averages will be computed). We further denote by $\Gamma_h^0:=\cup\{e: e\in\Ei\}$ the interior skeleton, by $\Gamma_h^b:=\cup\{e: e\in\Eb\}$ the boundary counterpart, and by $\Gamma_h:=\Gamma_h^0\cup\Gamma_h^b$ the full skeleton.

For $e \in \Ei$ we fix $\mu:=\mu_e$ to be one of the two unit normals to $e$ in $\Omega$; this choice is irrelevant for the discussion below. 
For $e \in \Eb$ we set $\mu$ to be the outward pointing unit normal to $e$.
Given $v_h\in \Vhk$, we denote
its piecewise gradient by $\nabla_h v_h$ and the {\it jumps} of $v_h$ and $\nabla_h v_h$
across any edge $e\in \calE_h$ by
\begin{equation}\label{jumps}
    [v_h]:= 
  \begin{cases}
  v_h^- -v_h^+ & e\in \Ei \\
  v_h^-             & e\in \Eb
  \end{cases},
  \qquad
  [\nabla_h v_h]:= 
  \begin{cases}
  \nabla_h v_h^- -\nabla_h v_h^+  & e\in \Ei \\
  \nabla_h v_h^-  & e\in \Eb 
  \end{cases},
\end{equation}
where $v_h^{\pm}(x) = \modif{\lim_{s \to 0^+}} v_h(x\pm s~\mu_e)$ for $x\in e$. 
The {\it averages} of $v_h\in\Vhk$ and $\nabla_h v_h$ across an
edge $e\in\calE_h$ are given by
\begin{equation}\label{averages}
\{ v_h \}:=
\begin{cases}
  \frac 1 2 (v_h^+ + v_h^-) & e\in \Ei \\
  v_h^- & e\in \Eb
\end{cases},
\quad
\{ \nabla_h v_h \}:=
\begin{cases}
  \frac 1 2 (\nabla_h v_h^+ + \nabla_h v_h^-) & e\in \Ei \\
  \nabla_h v_h^- & e\in \Eb
\end{cases}.
\end{equation}  
Motivated by the dG formulation of the bi-harmonic problem and \eqref{e:frob}, given two positive stabilization parameters $\gamma_0$ and $\gamma_1$,  we define the discrete energy $E_h^0$ on $\lbrack \Vhk \rbrack^3$  by \cite{BoNoNt}
\begin{equation} \label{E:SingleLayerDE}
\begin{aligned}
E_h^0[y_h] :& = \frac{1}{2}\|D_h^2y_h\|_{ \lom }^2 - ( f, y_h )_{ \lom }
\\
& -  (\{ \nd \nabla_h y_h \}, \left[ \nabla_h y_h \right] )_{\LtGi} + (\{ \nd \Delta_h y_h\}, \left[ y_h \right])_{\LtGi} \\
& -  (\{ \nd \nabla_h y_h \}, \nabla_h y_h - \Phi )_{\LtGb} + (\{ \nd \Delta_h y_h\}, y_h - g )_{\LtGb} \\
& + \frac{\gamma_1}{2} \|h^{-1/2} [\nabla_h y_h]\|_{\LtGi}^2 + \frac{\gamma_0}{2} \|h^{-3/2} [ y_h]\|_{\LtGi}^2 \\
&+ \frac{\gamma_1}{2} \|h^{-1/2} (\nabla_h y_h - \Phi) \|_{\LtGb}^2 + \frac{\gamma_0}{2} \|h^{-3/2} (y_h - g) \|_{\LtGb}^2 ,
\end{aligned}
\end{equation}
We show in \cite{BoNoNt} that, if $\gamma_0$, $\gamma_1$ are chosen sufficiently large and $\partial_D\Omega\ne\emptyset$, then the discrete energy $E_h^0$ is coercive with respect to the following dG quantity defined for $v_h \in \lbrack \Vhk \rbrack^3$
\begin{equation} \label{E:energynormN}
\begin{aligned}
||| v_h |||_{E}^2 :&=   \| D_h^2 v_h\|_{L^2(\Omega)}^2
+ \| h^{-1/2} \left[ \nabla_h v_h \right] \|_{\LtGi}^2 
+ \|h^{-3/2} [ v_h]\|_{\LtGi}^2  \\
&+ \| h^{-1/2} ( \nabla_h v_h - \Phi ) \|_{\LtGb}^2 
+ \|h^{-3/2} ( v_h - g ) \|_{\LtGb}^2,
\end{aligned}
\end{equation}
in the sense that $E_h^0[v_h] \leq c$ implies $|||v_h|||_E \leq C$ uniformly in $h$.
When $\partial_D \Omega = \emptyset$, the coercivity of $E_h^0$ holds on the subspace of $\Vhk$ consisting of functions with vanishing zero (mean value) and first moments.

Notice that the boundary conditions \eqref{E:bc} are enforced using the Nitsche method: if $E_h^0[y_h] \leq c$, then the coercivity property $|||v_h|||_E \leq C$ implies $\|y_h-g\|_{\LtGb}\leq C h^{3/2}$ and $\|\nabla_h y_h - \Phi\|_{\LtGb}\leq C h^{1/2}$,
whence \eqref{E:bc} is recovered as $h \to 0$. This allows us to avoid incorporating the boundary conditions explicitly on the discrete admissible set
\begin{equation}\label{E:admissible}
\Ahk := \left\{ 
\begin{aligned}
v_h \in [\Vhk]^3 : \ \sum_{T \in \mathcal T_h} \left| \int_T [\nabla_h v_h]^{\mathsf{T}}\nabla_h v_h-I \right| \leq \delta
\end{aligned}
\right\},
\end{equation}
where $\delta = \delta(h)$ is such that $\delta(h) \to 0$ as $h \to 0$.
By doing so, unlike for  the Kirchhoff approach \cite{BoNoNt}, technical issues related to compatibility between the isometry constraint \eqref{E:isometry} and boundary conditions \eqref{E:bc} are by-passed.
We also point out that imposing \eqref{E:isometry} at the vertices within the Kirchhoff approach is further relaxed to an elementwise average constraint in \eqref{E:admissible}, which is achievable with the gradient flow described below.

\medskip
{\bf $\Gamma$-convergence:}
We prove convergence of almost global minimizers in \cite[Section 5]{BoNoNt} provided
$\partial_D \Omega \not = \emptyset$. Let $y_h\in\Ahk$ be a sequence of almost global minimizers of $E_h$, i.e.
		$$
		E_h^0[y_h] \leq \inf_{w_h \in \Ahk} E_h^0[w_h] + \epsilon(h) \leq C,
		$$ 
		where $\epsilon(h) \rightarrow 0$, $\delta(h) \to 0$ as $h \rightarrow 0$ and 
		$C$ is a constant independent of $h$.
		Then, $\left\{ y_h \right\}_{h>0}$ is precompact in $[\lom]^3$ 
		and every cluster point $y$ of the sequence $\{y_h\}_{h>0}$
                satisfies
                $y \in [\htwom]^3, y=g, \nabla y = \Phi$ on $\partial_D \Omega$ and $y$ is an isometry and a global minimizer of $E^0$.
		Moreover, up to a subsequence (not relabeled), the discrete energy converges
		$$
		\lim_{h \rightarrow 0}E_h^0[y_h]=E^0[y].
		$$

\medskip
{\bf Discrete Euler-Lagrange equations:} In order to write discrete versions of \eqref{E:linear-isometry} and \eqref{E:EL}, we must first realize that the Dirichlet boundary conditions \eqref{E:bc} are imposed weakly via the Nitsche method. This changes the very notion of discrete tangent plane. We write the {\it discrete linearized isometry constraint}
\begin{equation}\label{e:isom_cond_local}
\modif{L_T[v_h;y_h] := \int_T  \Big( [\nabla v_h]^{\mathsf{T}} \nabla y_h + [\nabla y_h]^{\mathsf{T}} \nabla v_h \Big)= 0 \qquad\forall \, T\in\Th}
\end{equation}
for variations $v_h\in[\Vhk]^3$ of a minimizer $y_h\in\Ahk$.
This defines the {\it discrete tangent plane} $\calF_h[y_h]$ at $y_h$ to be
the set of functions $v_h\in[\Vhk]^3$ satisfying \eqref{e:isom_cond_local}
without boundary conditions. The discrete version of \eqref{E:EL} is $\delta E_h^0[y_h](v_h)=0$ for all $v_h\in\calF_h[y_h]$ and reads
\begin{equation}\label{E:discrete-EL}
  a_h(y_h,v_h) = (f,v_h)_{L^2(\Omega)} + \ell_h(v_h)
  \qquad\forall \, v_h\in\calF_h[y_h],
\end{equation}
where $a_h(\cdot,\cdot)$ is the bilinear form corresponding to \eqref{E:SingleLayerDE} and is given for $v_h,w_h \in [\Vhk]^3$ by
\begin{equation}\label{E:bilinear}
\begin{aligned}
  a_h\big(w_h,v_h\big) &:= \big(D_h^2 w_h, D_h^2 v_hv\big)_{ \lom }
  \\& - \big( \{ \nd \nabla_h w_h \}, \left[ \nabla_h v_h \right] \big)_{{L^2(\Gamma_h})}
	 - \big( \{ \nd \nabla_h v_h \}, \left[ \nabla_h w_h \right] \big)_{{L^2(\Gamma_h)}} \\
	 &+ \big( \{ \nd \Delta_h w_h \}, \left[ v_h \right] \big)_{{L^2(\Gamma_h)}}
	 + \big( \{ \nd \Delta_h v_h \}, \left[  w_h \right] \big)_{{L^2(\Gamma_h)}} \\
	 & + \gamma_1 \big(h^{-1} \left[ \nabla_h w_h \right], \left[ \nabla_h v_h \right] \big)_{{L^2(\Gamma_h)}} + \gamma_0 \big(h^{-3} \left[ w_h \right], \left[ v_h \right] \big)_{{L^2(\Gamma_h)}},
\end{aligned}
\end{equation}
and $\ell_h$ is the linear form that enforces the boundary conditions in the
Nitsche's sense
\begin{align*}
  \ell_h(v_h) &:= -\big( \nd \nabla_h v_h, \Phi  \big)_{\LtGb}  + \big( \nd \Delta_h v_h ,   g  \big)_{\LtGb}
  \\ &+ \gamma_1 \big(h^{-1} \Phi , \nabla_h v_h  \big)_{\LtGb} + \gamma_0 \big( h^{-3} g ,  v_h  \big)_{\LtGb}.
\end{align*}
To see this, simply note that we use the full skeleton $\Gamma_h$ in \eqref{E:bilinear}
and combine the four terms in $\ell_h(v_h)$ with similar terms on the boundary skeleton
$\Gamma_h^b$ in \eqref{E:bilinear} to arrive at the following equivalent form of \eqref{E:discrete-EL} where the Nitsche's approach is apparent:
\begin{equation}\label{E:Nitsche}
\begin{aligned}
  \big(D_h^2 y_h, D_h^2 v_h\big)_{ \lom } &-
  \big( \{ \nd \nabla_h y_h \}, \left[ \nabla_h v_h \right] \big)_{{L^2(\Gamma_h})}
	 - \big( \{ \nd \nabla_h v_h \}, \left[ \nabla_h y_h \right] \big)_{{L^2(\Gamma_h^0)}} \\
	 &+ \big( \{ \nd \Delta_h y_h \}, \left[ v_h \right] \big)_{{L^2(\Gamma_h)}}
	 + \big( \{ \nd \Delta_h v_h \}, \left[  y_h \right] \big)_{{L^2(\Gamma_h^0)}} \\
	 & + \gamma_1 \big(h^{-1} \left[ \nabla_h y_h \right], \left[ \nabla_h v_h \right] v\big)_{{L^2(\Gamma_h^0)}} + \gamma_0 \big(h^{-3} \left[ y_h \right], \left[ v_h \right] \big)_{{L^2(\Gamma_h^0)}}
         \\
         & 
         -\big( \nd \nabla_h v_h, \nabla_h y_h-\Phi  \big)_{\LtGb}  + \big( \nd \Delta_h v_h ,  y_h- g  \big)_{\LtGb}
         \\& + \gamma_1 \big( h^{-1}(\nabla_h y_h - \Phi) , \nabla_h v_h  \big)_{\LtGb} + \gamma_0 \big( h^{-3}(y_h - g) ,  v_h  \big)_{\LtGb}
         \\& = \big(f,v_h\big)_{\lom}.
\end{aligned}
\end{equation}
  
\medskip
{\bf Discrete gradient flow:}
In order to construct discrete minimizers $y_h\in[\Vhk]^3$ of \eqref{E:SingleLayerDE} satisfying \eqref{e:isom_cond_local}, or solutions of \eqref{E:discrete-EL}, we employ the following discrete relaxation dynamics with pseudo-time step $\tau>0$. Given a current deformation $y_h^n \in \lbrack \Vhk \rbrack^3$ at iteration $n \in \mathbb N$, we seek the new iterate $y_n^{n+1}:=y_h^n + \delta y_n^{n+1}$ with correction $\delta y_h^{n+1} \in \calF_h[y_h^n]$ satisfying
\begin{equation}\label{E:grad-flow1}
  \begin{aligned}
\tau^{-1} \big(\delta y_h^{n+1},v_h \big)_{H_h^2}
+ \ a_h \big(\delta y_h^{n+1},v_h \big) &=  - a_h \big( y_h^{n},v_h \big)
\\ & + \ (f, v_h)_{\lom} + \ \ell_h(v_h)
\quad\forall \, v_h\in\calF_h[y_h^n].
  \end{aligned}
\end{equation} 
This is a discrete version of \eqref{E:gf-1} in $[\Vhk]^3$ with variations $\delta y_h^{n+1}$ tangent
to $y_h^n$ and underlying metric induced by the discrete $\htwom$-inner product
\modif{$\big(\cdot,\cdot\big)_{H_h^2}$} and corresponding norm $|||\cdot|||_{H_h^2}$,
where
\begin{equation}\label{E:inner-product}
\begin{aligned}
  \big(v_h,w_h\big)_{H_h^2} &:= \big(D_h^2v_h, D_h^2w_h\big)_\lom + \varepsilon   \big(v_h,w_h\big)_\lom
  \\ & + \big(h^{-1}[\nabla v_h], [\nabla w_h] \big)_{L^2(\Gamma_h^0)} + \big(h^{-3}[v_h],[w_h] \big)_{L^2(\Gamma_h^0)}.
\end{aligned}
\end{equation}
Notice that the presence of the $\varepsilon$-term as in \eqref{e:H2} ensures that $\big(\cdot,\cdot\big)_{H_h^2}$ is indeed an inner product in $[\Vhk]^3$ and that therefore \eqref{E:grad-flow1} has a unique solution even when $\partial_D \Omega = \emptyset$. 

An important property of gradient flows is that the resulting deformations $y_h^{n+1} = y_h^n + \delta y_h^{n+1} \in \lbrack \Vhk \rbrack^3$ decrease the discrete energy strictly
provided $\delta y_h^{n+1} \ne 0$ \cite[Lemma 3.2]{BoNoNt}
\begin{equation}\label{E:discrete-energy-decrease}
 \frac{1}{\tau} |||\delta y_h^{n+1} |||_{H^2_h}^2 + E_h^0[y_h^{n+1}] \le E_h^0[y_h^n];
\end{equation}
this is the discrete counterpart of \eqref{E:energy-decrease}.
This also shows that the sequence $y_h^n$ converges to a local minimizer
$y_h\in[\Vhk]^3$ of $E_h^0$.
In addition, if \modif{$\delta \geq \big(1 + c_1E_h^0[y_h^0] + c_2 R(g,\Phi,f)\big)\tau$ with $R(g,\Phi,f)=\|g\|_{H^1(\Omega)}^2+ \|\Phi\|_{H^1(\Omega)}^2 + \|f\|_{L^2(\Omega)}^2$,} then the linearized isometry constraint \eqref{e:isom_cond_local}, together with \eqref{E:discrete-energy-decrease} and Friedrichs inequality, guarantees $y_h^n \in \Ahk$ for all $n\ge1$; hence, $y_h\in \Ahk$. We refer to \modif{\cite[Lemma 3.2]{BoNoNt}} for additional details.

\section{Simplified Bilayer Model and its dG Approximation} \label{S:SimpleBilayer}

We now briefly recall the simplified bilayer model from \cite{BaBoMuNo}, and adjust the dG method of \cite{BoNoNt} to this model.

\medskip
{\bf Continuous energy:} Since $y$ is a pointwise isometry, i.e. $y$ satisfies \eqref{E:isometry}, in view of \eqref{e:frob} the energy functional \eqref{E:bilayer_energy} for the bilayer plate can be rewritten as
\begin{equation}\label{E:bilayer_energy2}
  E^1[y] := \frac{1}{2} \int_{\Omega} |D^2 y|^2 - \int_\Omega H : Z + \frac12\int_\Omega |Z|^2 - \int_{\Omega} f y,
\end{equation}
where the matrix function $Z\in\mathbb{R}^{2\times2}$ is referred to as a \emph{spontaneous curvature} and encodes the mismatch between the two constituent materials of the bilayer plate. 
Motivated by the applications presented in Section~\ref{S:NumericalExperiments}, we focus on the case $f=0$ in the discussion below.
Moreover, using the expression \eqref{E:normal} for the unit normal, valid for isometries, the second fundamental form reads $H=(\partial_{ij} y \cdot (\partial_1 y \times \partial_2 y))_{ij}$ and leads to
\begin{equation} \label{E:bilayer_energy3}
E^1[y] = \frac{1}{2} \int_\o |D^2y|^2  - \sum_{ij} \int_\Omega \partial_{ij} y  \cdot (\partial_1 y \times \partial_2 y) \,  z_{ij}  + \frac12\int_\Omega |Z|^2.
\end{equation}

The first variation $\delta E^1[y](v)$ of $E$ at $y \in [H^2(\Omega)]^3$ in the direction
$v \in \calF[y]$ reads
\begin{align*}
  \delta E^1[y](v) := \int_{\Omega} D^2y: D^2v &- \sum_{ij} \int_{\Omega} \partial_{ij} v  \cdot ( \partial_1 y \times \partial_2 y) \, z_{ij}
  \\
  & - \sum_{ij} \int_{\Omega} \partial_{ij} y  \cdot ( \partial_1 y \times \partial_2 v + \partial_1 v \times \partial_2 y ) \, z_{ij},
\end{align*}
where we recall that $\calF[y]$ is the set of functions $v\in[H^2(\Omega)]^3$ with vanishing Dirichlet boundary conditions and satisfying the linearized isometry constraint \eqref{E:linear-isometry}. To obtain an equivalent expression better suited to numerical approximations, we combine \eqref{E:hessian} with the vector identity \modif{$(a \times b) \cdot (c \times d) = (a\cdot c) (b \cdot d) - (a \cdot d)(b \cdot c)$} to realize that
\begin{align*}
  \partial_{ij} y & \cdot ( \partial_1 y \times \partial_2 v) =
   h_{ij} \ ( \partial_1 y \times \partial_2 y) \cdot  (\partial_1 y \times \partial_2 v) = 0
  \\
  \partial_{ij} y & \cdot (\partial_1 v \times \partial_2 y ) =
   h_{ij} \ ( \partial_1 y \times \partial_2 y) \cdot  (\partial_1 v \times \partial_2 y) = 0
\end{align*}
because $\partial_1 y \cdot \partial_2 y=0$ and $\partial_iy\cdot\partial_i v=0$ according to \eqref{E:isometry} and \eqref{E:linear-isometry} respectively.
Therefore, the expression of  $\delta E^1[y](v)=0$ simplifies to
\begin{equation} \label{E:bilayer_simple}
	\delta E^1[y](v) = \int_{\o} D^2y: D^2v - \sum_{ij} \int_{\o} \partial_{ij} v  \cdot  ( \partial_1 y \times \partial_2 y) \, z_{ij} = 0,
\end{equation}
whenever $y\in [H^2(\Omega)]^3$ is an isometry that minimizes \eqref{E:bilayer_energy3} and $v\in\calF[y]$. In analogy with \eqref{E:gf-1}, we can now find $y$ as the asymptotic limit of the solution of the $H^2$-gradient flow
\begin{equation}\label{E:gf-2}
  \big( \partial_t y, v \big)_{H^2(\Omega)} + \delta E^1[y](v) = 0
  \quad\forall \, v\in \calF[y]
  \qquad\Rightarrow\qquad \frac{d}{dt} E^1[y]\le 0.
\end{equation}
We observe that the first term in \eqref{E:bilayer_simple} dominates the second one and is already present in the first variation \eqref{E:EL} of the single layer energy functional $E^0$. We exploit this next.

\medskip
{\bf Discrete energy:}
In order to obtain a simple yet efficient discretization of \eqref{E:bilayer_simple}, we take advantage of the good properties of the dG discretization of \Cref{S:SingleLayer}. We thus discretize the first term in \eqref{E:bilayer_simple} according to \eqref{E:discrete-EL} and the second one elementwise to arrive at
\begin{equation}\label{E:bilayer-EL}
  a_h(y_h,v_h) = \ell_h(v_h) +
   \sum_{ij} \sum_{ T \in \mathcal T_h}  \int_T  \partial_{ij} v_h  \cdot ( \partial_1 y_h \times \partial_2 y_h) \, z_{ij} 
  \qquad\forall \, v_h\in\calF_h[y_h].
\end{equation}
We emphasize that this {\it nonlinear} discrete scheme entails only a piecewise computation of the additional term without taking into account any possible jumps of $\partial_{ij} v_h$ and of $\partial_1 y_h \times \partial_2 y_h$. They are indeed already incorporated into the bilinear form $a_h(\cdot,\cdot)$ defined in \eqref{E:bilinear} and provide a good approximation of the Hessian $D^2y$ \cite{BoNoNt}. The discrete energy $E_h^1$ associated with \eqref{E:bilayer-EL} reads
\[
E_h^1[y_h] := E_h^0[y_h] - \sum_{ij} \sum_{ T \in \mathcal T_h}  \int_T  \partial_{ij} y_h  \cdot ( \partial_1 y_h \times \partial_2 y_h) \, z_{ij}.
\]
We point out that \eqref{E:bilayer-EL} is not the Euler-Lagrange equation of $E_h^1$ because the orthogonality conditions leading to \eqref{E:bilayer_simple} may not be valid for $y_h \in\Ahk$. We deal with \eqref{E:bilayer-EL} below.

\medskip
{\bf Discrete gradient flow:} To compute a solution $y_h\in\Ahk$ of \eqref{E:bilayer-EL} we propose a discrete version of \eqref{E:gf-2} with $\delta E^1_h[y_h^{n+1}](v_h)$ replaced by \eqref{E:bilayer-EL} and its rightmost term treated explicitly; compare to \eqref{E:grad-flow1}. 
Given $y_h^n \in \Ahk$, we thus seek $\delta y_h^{n+1} \in \calF[y_h^n]$ such that
\begin{equation}\label{E:grad-flow2}
\begin{aligned}
\tau^{-1}\big(\delta y_h^{n+1},v_h \big)_{H_h^2}
& + \ a_h \big(\delta y_h^{n+1},v_h \big) = - \ a_h \big( y_h^{n},v_h \big) \ + \ \ell_h(v_h) \\
& + \sum_{ij} \sum_{ T \in \mathcal T_h}  \int_T z_{ij} \ \partial_{ij}  v_h  \cdot  ( \partial_1 y_h^n \times \partial_2 y_h^n)\qquad\forall \ v_h\in\calF[y_h^n],
\end{aligned}
\end{equation}
and set $y_h^{n+1} := y_h^n + \delta y_h^{n+1}$. This {\it linear} algorithm is used for the simulations in \Cref{S:NumericalExperiments}.

Several comments are in order. It is not clear that \eqref{E:grad-flow2} reduces the elastic energy $E_h^1[y_h^n]$ and $y_h^{n+1}\in\Ahk$. To show these crucial properties we need to quantify for $y_h^n$ the lack of orthogonality leading to \eqref{E:bilayer_simple}. Moreover, we have to quantify the effect of the explicit treatment of the last term in \eqref{E:grad-flow2}.
Finally, the isometry constraint \eqref{E:isometry} is not valid for $y_h^n\in\Ahk$, for which we allow an isometry defect. In view of \eqref{e:isom_cond_local}, we deduce
\[
\int_T [\nabla y_h^{n+1}]^{\mathsf{T}} \nabla y_h^{n+1} =
\int_T [\nabla y_h^n]^{\mathsf{T}} \nabla y_h^n +
\int_T [\nabla \delta y_h^{n+1}]^{\mathsf{T}} \nabla \delta y_h^{n+1}
\ge \int_T [\nabla y_h^n]^{\mathsf{T}} \nabla y_h^n
\quad\forall \, T\in\Th,
\]
whence exploiting telescopic cancellation we obtain
\[
|T|^{-1}\int_T [\nabla y_h^n]^{\mathsf{T}} \nabla y_h^n \ge |T|^{-1}\int_T [\nabla y_h^0]^{\mathsf{T}} \nabla y_h^0 = I
\]
provided $y_h^0$ is an isometry in the sense that the last equality holds. This implies that the average of $\partial_i y_h^n \ge 1$ over each $T\in\Th$ and \eqref{E:grad-flow2} is well-defined. However, the vector $\partial_1 y_h^n \times \partial_2 y_h^n$ may not have unit norm. These observations were instrumental in \cite{BaBoNo} to redefine the normal vector as $\frac{\partial_1 y_h^n}{|\partial_1 y_h^n|} \times \frac{\partial_2 y_h^n}{|\partial_2 y_h^n|}$ and prove $\Gamma$-convergence of the discrete energy. We anticipate that a rigorous study of the dG approach \eqref{E:grad-flow2} to the bilayer model may need to take advantage of the \modif{properties of the discrete Hessian $H_h[y_h]$} described in \cite{BoNoNt}, as well as the \modif{quasi-orthogonality} relations
\[
\sum_{T\in\Th} \Big| \int_T \partial_1 y_h^n \cdot \partial_2 y_h^n \, \Big| \le \delta
\]
and
\[
\int_T \partial_i y_h^n\cdot\partial_i v_h = 0
\qquad\forall \, v_h\in\calF_h[y_h^n], \quad\forall \, T\in\Th,
\]
ensuing from \eqref{E:admissible} and \eqref{e:isom_cond_local} respectively.
We do not investigate here any theoretical properties of \eqref{E:grad-flow2} but rather explore its performance on several numerical experiments, some from \cite{BaBoNo, BaBoMuNo} and some new. This is carried out in the next section.

\section{Numerical Experiments}\label{S:NumericalExperiments}

In this section we explore several examples motivated by the work in \cite{BaBoNo,BaBoMuNo} and experimental work in \cite{AlBSm,Simpson,JHZ} in order to further understand the computational performance of the dG method \eqref{E:grad-flow2}, which extends our method from \cite{BoNoNt} to bilayer plates. We aim to verify whether the simple extension of \Cref{S:SimpleBilayer} can lead to relevant simulations that capture the essential physical properties of the bilayer bending problem seen in lab experiments. To this end, we challenge our algorithm in a variety of settings exploring the effect of the spontaneous curvature matrix $Z$, the boundary conditions and the midplate aspect ratio. All these factors are crucial for engineering applications as they can be appropriately combined to achieve desirable equilibrium configurations under suitable thermal, electrical, or chemical actuation. We emphasize that the model discussed here does not assume small deformations and therefore copes with geometrically nonlinear deformations. 

\subsection{Boundary Conditions}\label{ss:BC}
Before we proceed with the specific examples, it is worth discussing briefly our approach for boundary conditions (B.C.). Our experiments fall into two main categories:

\begin{enumerate}[$\bullet$]
\item \textbf{Dirichlet B.C. on $\partial_D \Omega$: Case $\varepsilon=0$.} We impose the Dirichlet condition \eqref{E:bc} via a Nitsche approach in \eqref{E:Nitsche}, whose left-hand side is the variational derivative of $E_h^0[y_h]$ in \eqref{E:SingleLayerDE}. In fact, the left-hand side of \eqref{E:Nitsche} contains boundary terms on $\Gamma_h^b$ with the quantities $\nabla_h y_h - \Phi$ and $y_h-g$ but, incidentally, does not contain terms on $\partial \Omega \setminus \Gamma_h^b$. Therefore, the discrete space $[\Vhk]^3$ does not include Dirichlet boundary conditions. The corresponding gradient flow \eqref{E:grad-flow2} is linear and coercive thanks to a Friedrichs inequality \cite{Brenner},\cite[Corollary 2.2]{BoNoNt} even when $\varepsilon = 0$ in \eqref{E:inner-product}; hence \eqref{E:grad-flow2} admits a unique solution. Moreover, our $\Gamma-$convergence analysis of  \cite{BoNoNt} guarantees that the limiting deformation $y$ satisfies \eqref{E:bc} on $\partial_D \Omega$ for single-layer plates. All the numerical simulations with $\partial_D \Omega\ne\emptyset$ are therefore performed with $\varepsilon = 0$.

\smallskip          
\item \textbf{Free B.C.: Case $\varepsilon>0$.}  When $\partial_D \Omega = \emptyset$, we allow the thin plate to deform without any boundary restrictions. We thus realize that the kernel of the bilinear form on the left-hand side of \eqref{E:grad-flow2} is non-trivial unless $\varepsilon > 0$. The characterization of such kernel with $\varepsilon=0$ and 
linearized isometry constraint is not clear. Nevertheless, to guarantee that the resulting linear system is uniquely solvable, we add a zero order $\varepsilon$-term to the discrete semi-inner product $(\cdot,\cdot)_{H_h^2}$ defined in \eqref{E:inner-product}. Stationary configurations of \eqref{E:grad-flow2} are unaffected by the $\varepsilon$-term. In practice, we take $\varepsilon=10^{-2}$ in our examples below.
\end{enumerate}

\subsection{Implementation of the Gradient Flow}

We now briefly describe some implementation aspects of the gradient flow. \modif{Except for the convergence analysis of Section \ref{S:clamped}, all subsequent experiments are performed on 5 uniform refinements of the plate, resulting in 1024 cells, and with a pseudo-time step $\tau = 5 \cdot 10^{-3}$. The associated finite element spaces are given by \eqref{E:discrete-space} for the case of quadrilaterals and $k=2$, i.e. using $\mathbb Q^2$ polynomials in the reference element $\widehat T$.}
 Based on parameter studies performed for the single layer model in \cite{BoNoNt},  we use rather large values for the stabilization parameters: $\gamma_0 = 5 \cdot 10^3$ and $\gamma_1 = 1.1 \cdot 10^3$. This choice is not exclusively dictated by stability considerations, as is customary for interior penalty dG \cite{BoNo,Pryer,Riv}, but primarily by a balance between the discrete initial energy $E_h^0[y_h^0]$ and the fictitious time-step $\tau$ of the discrete gradient flow, which determines the isometry defect $\delta$ in \eqref{E:admissible} according to \modif{$\delta \ge \big(1+c\big)\tau$ for the single layer model, where $c$ is a constant depending on the initial energy $E_h^0[y_h^0]$, the data $f,g,\Phi$ and the domains $\Omega, \partial_D \Omega$.} It is important to notice that the magnitudes of $\gamma_0,\gamma_1$ affect the weak imposition of Dirichlet conditions alla Nitsche and the value of $E_h^0[y_h^0]$, thereby making the choice of $\gamma_0,\gamma_1$ a \modif{critical} aspect of our method. Different choices may influence the deformation flow and potentially lead to cases where local discrete minimizers are attained instead of global ones. However we recall that for the single layer system, once $\gamma_0,\gamma_1$ are fixed, $\Gamma$-convergence guarantees that discrete global minimizers converge towards exact global minimizers.

The stopping criteria for the gradient flow is
\begin{equation}\label{e:tolerance_algo}
	\big| E_h^1[y_h^{n+1}]-E_h^1[y_h^n] \big| < 10^{-6} \ \tau = 5 \cdot 10^{-9}.
\end{equation}
We declare such a deformation $y_h^{n+1}$ to be our equilibrium deformation.

Lastly, to implement the linearized isometry constraint \eqref{e:isom_cond_local}, we use a piecewise constant symmetric Lagrange multiplier matrix $\lambda_h^{n+1}\in[\Vh^0]^{2\times2}$ with $3$ components corresponding to the distinct elements of the
symmetric bilinear form $L_T[v_h; y_h^n]\in[\Vh^0]^{2\times2}$ of \eqref{e:isom_cond_local} \modif{for $T\in\Th$. If $\{ \varphi_i \}$ stands for the standard Lagrange basis of $[\mathbb V_h^k]^3$ (i.e. piecewise $\mathbb P^k$ or $\mathbb Q^k$ in the reference element) and $\{ \Psi_i \}$ for the basis of $[\mathbb V_h^0]^{2\times 2}$ consisting of piecewise constant $2 \times 2$ matrices, then the bilinear form associated with $\lambda_h^{n+1}$ is
\[
b^n_T(\lambda_h^{n+1},\varphi_i):=\int_T \lambda_h^{n+1} : \Big( [\nabla \varphi_i]^{\mathsf{T}} \nabla y_h^n +
    [\nabla y_h^n]^{\mathsf{T}} \nabla \varphi_i  \Big)
    \qquad\forall \, T\in\Th.
\]
Moreover, we denote by $\Lambda_h^{n+1}$ the vector representation of $\lambda_h^{n+1}$ in the basis $\{ \Psi_i \}$ and by $Y^{n+1}$ the vector representation of $y_h^n$ in the basis $\{ \varphi_i\}$; hence $\delta Y^{n+1} := Y^{n+1}-Y^n$. Therefore, the augmented linear system corresponding to \eqref{E:grad-flow2} with Lagrange multiplier is a saddle point system and reads
\begin{equation}\label{e:saddle_explicit}
\begin{pmatrix}
\tau^{-1} M + A & (B^n)^{\mathsf{T}} \\
B^n & 0 
\end{pmatrix}
\begin{pmatrix}
\delta Y^{n+1} \\
\Lambda^{n+1} 
\end{pmatrix}
=
\begin{pmatrix}
-A Y^n + F^n +G\\
0
\end{pmatrix},
\end{equation}
where the coefficients of $M, A$ and $B^n$ are
\begin{equation}\label{e:M}
M_{ij} := (\varphi_j,\varphi_i)_{H^2_h}, \quad A_{ij}:= a_h(\varphi_j,\varphi_i), \quad B^n_{ij} := \sum_{T\in \mathcal T_h} b^n_T(\Psi_i,\varphi_j)
\end{equation}
whereas those of the right-hand side of \eqref{e:saddle_explicit} are 
\begin{equation}\label{e:RHS}
  F^n_i :=  \sum_{k,l}\sum_{T\in \mathcal T_h} \int_T z_{kl} \, \partial_{kl} \varphi_i \cdot (\partial_1 y^n_h \times \partial_2 y_h^n), \qquad G_i := \ell_h(\varphi_i).
\end{equation}
}

\begin{algorithm}[h!]
\modif{
\caption{Gradient Flow} \label{algo:insertion-sort}
\begin{algorithmic}[1]
\State \emph{// Input} 
\State \qquad $y_h^0 \leftrightarrow Y^0$ \Comment{Initial deformation function $\leftrightarrow$ vector} 
\State \qquad $\gamma_0$, $\gamma_1$ \Comment{Jump penalization parameters}
\State \qquad $\tau$ \Comment{Pseudo-time step} 
\algrule
\State \emph{// Initialization}
\State $n=0$
\State Assemble $M$, $A$, $B^0$ according to \eqref{e:M}
\State Assemble $F^0$, $G$ according to \eqref{e:RHS}
\State Assemble $\mathcal A = \tau^{-1} M +A$ and compute its LU decomposition
\State Compute $E_h^1[y_h^0]$
\State \emph{// Main Loop}
\Repeat 
\State Solve for $\Lambda^{n+1}$ according to \eqref{e:shur_B} 
\State Solve for $\delta Y^{n+1}$ according to \eqref{e:shur_Y}
\State Compute $E_h^1[y_h^{n+1}]$
\State $n \gets n+1$
\Until{$\abs{E_h^1[y_h^{n}]-E_h^1[y_h^{n-1}]} <  10^{-6} \tau$}
\end{algorithmic}
}
\end{algorithm}

\modif{
We propose a Schur complement iterative method to solve \eqref{e:saddle_explicit}, namely
\begin{equation}\label{e:shur_B}
B^n \mathcal A^{-1} (B^n)^{\mathsf{T}} \Lambda^{n+1} = B^n \mathcal A^{-1} \big(-A Y^n + F^n +G \big), 
\end{equation}
and
\begin{equation}\label{e:shur_Y}
\delta Y^{n+1} = \mathcal A^{-1} \big(-A Y^n + F^n +G - (B^n)^{\mathsf{T}} \Lambda^{n+1} \big),
\end{equation}
where the matrix $\mathcal A:=\tau^{-1} M +A$ needs to be assembled only once and so the application of its inverse can be computed using a LU decomposition (computed once as well).
This results in an efficient inner solver used in conjunction with a conjugate gradient algorithm to compute $\Lambda^{n+1}$ in \eqref{e:shur_B}; efficient preconditioning of \eqref{e:shur_B} is still an open issue.
A pseudo-code for the full algorithm is given in Algorithm~\ref{algo:insertion-sort}.
Additional details can be found in \cite{BoNoNt}. The implementation of the dG method and simulations below have been carried out within the FEM software platform {\it deal.ii} \cite{dealii85,dealiipaper}. \modif{The resulting deformations are visualized with {\it Paraview} \cite{paraview}.}
}
 
\subsection{Clamped Plate: $Z=I$ and Comparison with Existing Methods}\label{S:clamped}
\begin{figure}[htb]
  \begin{center}
	\includegraphics[scale=0.45]{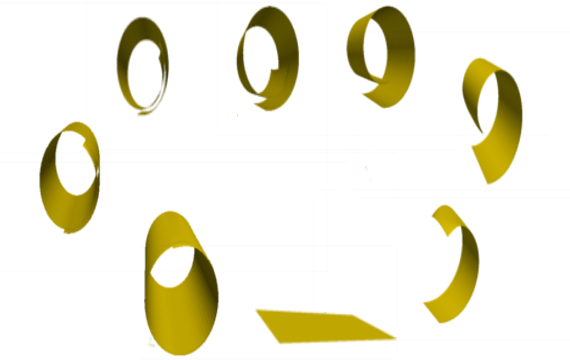}
  \end{center}
  \caption{\small Pseudo-evolution (counter-clockwise) towards the equilibrium of a clamped rectangular plate with spontaneous curvature $Z=I$. The bilayer plate is depicted at times 0.0, 0.5, 1.5, 3.0, 5.0, 8.0, 10.0, 20.0 $\times 10^{3}$ of the gradient flow. The plate reaches a cylindrical shape asymptotically, which is an absolute minimizer, without self-crossing.}
	\label{F:Identity}
\end{figure}

We start with the simplest case of a rectangular plate \modif{$\Omega = (-5,5) \times (-2,2)$,} clamped on the side $\left\{ -5\right\} \times [-2,2]$, with spontaneous curvature given by $Z=I$. \modif{The deformation with minimal energy corresponds to a cylinder of radius $1$ and energy 20 \cite{Schmidt}}. We report relevant iterations of the gradient flow in Figure \ref{F:Identity}. A cylindrical equilibrium configuration is reached confirming the results in \cite{BaBoNo,Schmidt}. It is worth mentioning that, in contrast to \cite{BaBoNo}, this final stage is attained now without any self-crossing of the plate, which is generally speaking difficult to avoid for a relaxation dynamics such as \eqref{E:grad-flow2}.

\FloatBarrier
\FloatBarrier

\begin{table}
\modif{
\begin{center}
\begin{tabular}{c||c|c||c|c||c|c}
Mesh  & \multicolumn{2}{c||}{\#4} & \multicolumn{2}{c||}{\#5}  & \multicolumn{2}{c}{\#6}  \\
\hline
\hline
Method & dG & K & dG & K & dG & K \\
Energy & 18.514 & 15.961 & 18.679 & 16.544 & 18.891 & No Convergence
\end{tabular}
\caption{Final energies using the Kirchhoff method (K) of \cite{BaBoNo} and the proposed dG method on $4$, $5$ and $6$ consecutive uniform refinements of the plate $(-5,5)\times (-2,2)$. The exact equilibrium energy is $20$ and corresponds to a cylindrical shape. We observe that dG method is at least $10\%$ more accurate. Moreover, unlike the dG method, the Kirchhoff method is not able to reach a stationary state for refinement $6$.}
\label{t:comp}
\end{center}
}
\end{table}

\modif{
We exploit this example to illustrate the improved accuracy and geometric flexibility of the dG method relative to the Kirchhoff approach in \cite{Bartels} for $Z \not  = 0$. 
A complete numerical study of this property when $Z=0$ is already presented in \cite[ Section 6]{BoNoNt}. We consider a sequence of $4,5$ and $6$ consecutive refinements of the plate $\Omega$, thus resulting in $256$, $1024$ and $4096$ cells with a total of
$7680$, $30720$ and $122880$ degrees of freedom; they are denoted by mesh \#4, \#5 and \#6 respectively. The pseudo-time step is set to $\tau  = 5 \cdot 10^{-3}$; we refer again to \cite{BoNoNt} for a convergence analysis with respect to the time-step. Table~\ref{t:comp} displays the final energies obtained with the proposed dG method and the Kirchhoff method of \cite[Table 6.1]{BaBoNo}. We observe that the dG method gives equilibrium deformations with energy $10\%$ more accurate that the Kirchhoff method. 
}

\subsection{Clamped Plate: Anisotropic Curvature}\label{S:anisotropy}

\begin{figure}[htb]
	\begin{center}
		\includegraphics[scale=0.42]{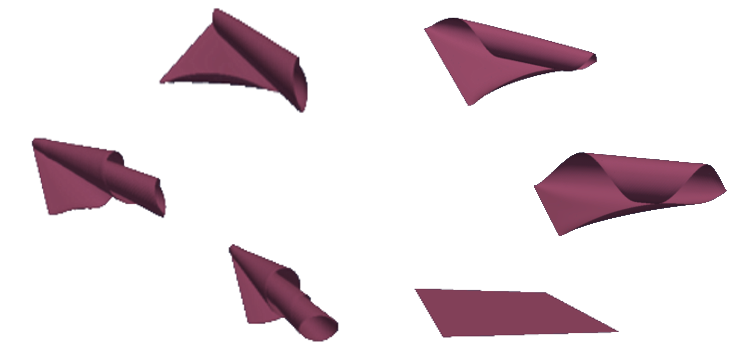}
	\end{center}
	\caption{\small Pseudo-evolution (counter-clockwise) towards the equilibrium of a clamped rectangular plate with principal directions forming an angle of $\pi/4$ with the coordinate axes and principal curvatures $5$ and $1$ (anisotropic spontaneous curvature). The bilayer plate is depicted at times 0.0, 0.3, 1.0, 10.0, 50.0, 170.0 $\times 10^{3}$ of the gradient flow and rolls to a conic shape.}
	\label{F:Anisotropic}
\end{figure}

We now explore the effect of anisotropic spontaneous curvature
$$
  Z =
	\left[
	\begin{matrix}
	3 & -2 \\
	-2 & 3
	\end{matrix}
	\right].
$$
        The plate \modif{$\Omega=(-2,2)\times(-3,3)$} is clamped along the side $\partial_D \Omega=[-2,2]\times \left\{-3 \right\}$. In contrast to Example \ref{S:clamped}, this choice of $Z$ corresponds to principal curvatures $5$ and $1$ and principal directions forming an angle $\pi/4$ with the coordinate axes. This, in conjunction with clamped boundary conditions, yields a plate that gradually ``rolls" into a conic shape rather than a cylindrical shape. We illustrate this deformation in \Cref{F:Anisotropic}.

\FloatBarrier
\FloatBarrier        

\subsection{Clamped Plate: Principal Curvatures of Opposite sign}

\begin{figure}[htb]
  \begin{center}
    \includegraphics[height=0.25\textheight, width=0.7\textwidth]{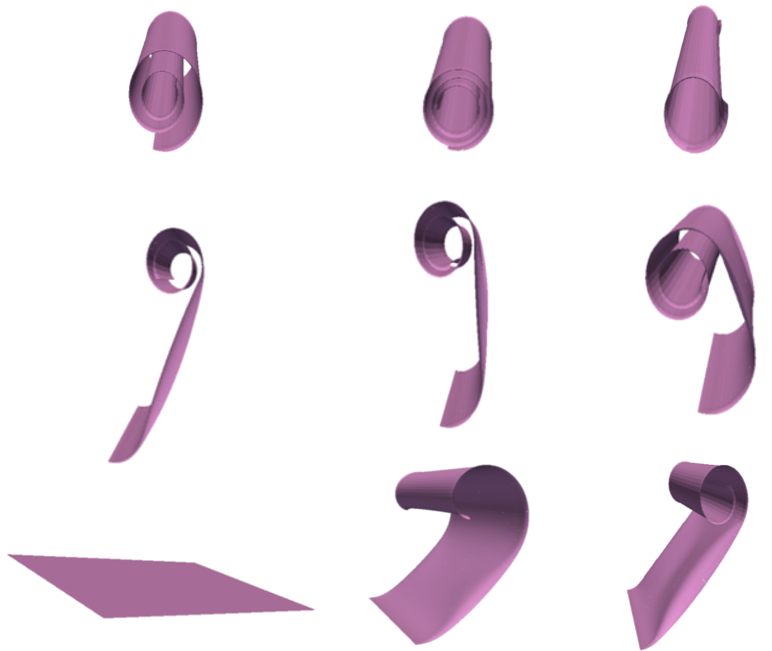}
  \end{center}
	\caption[Bilayer: Clamped Opposite Signs]{\small Pseudo-evolution (bottom-left to top-right) towards the equilibrium of a clamped rectangular plate with principal curvatures of opposite sign. The bilayer plate is depicted at times 0.0, 0.7, 2.5, 7.0, 10.0, 12.0, 14.0, 19.0, 67.0 $\times 10^{3}$ of the gradient flow. The equilibrium shape is a cylinder.}
	\label{F:Alternate}
\end{figure}

We consider the same plate as in Example \ref{S:anisotropy}, namely \modif{$\Omega=(-2,2)\times(-3,3)$} clamped along the side $\partial_D \Omega=[-2,2]\times \left\{-3 \right\}$, but with spontaneous curvature
$$
Z =
\left[
\begin{matrix}
-5 & 0 \\
0 & 5
\end{matrix}
\right].
$$
The principal curvatures have now opposite signs and the principal directions are the coordinate axes. At first the plate tries to bend in each coordinate direction according to the corresponding curvature sign. Eventually the longer side dominates and the structure attains the cylindrical configuration depicted in Figure\nobreakspace \ref {F:Alternate}.
\FloatBarrier
\FloatBarrier

\FloatBarrier
\begin{figure}[htb]
  \begin{center}
    \includegraphics[height=0.25\textheight, width=0.7\textwidth]{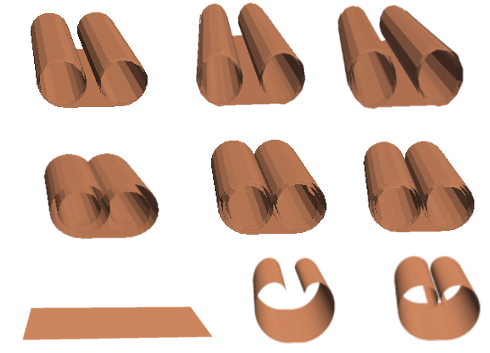}
   \end{center}
	\caption[Bilayer: Middle Clamped Anisotropic]{Pseudo-evolution (bottom-left to top-right) towards equilibrium of a rectangular plate clamped in the middle, with anisotropic spontaneous curvature. The bilayer plate is depicted at times 0.0, 0.2, 0.6, 6.0, 10.0, 12.0, 20.0, 60.0, 155.0 $\times 10^{3}$ of the gradient flow and self-intersects. The equilibrium configuration is two identical but disjoint cylinders.}
	\label{F:Anisotropic_Middle}
\end{figure} 
\FloatBarrier

\subsection{Plate Clamped in the Middle: Anisotropic Curvature}\label{S:middle}

We now explore an example that was motivated by \cite{AlBSm}. The rectangular plate is \modif{$\Omega = (-5,5)\times(-2,2)$} and we set the Dirichlet condition to be $g(x)=x$ and $\Phi(x)=0$ on the middle line \modif{$\left\{0 \right\}\times [-2,2]$ of $\Omega$.}
The plate is endowed with an anisotropic spontaneous curvature
$$
Z =
\left[
\begin{matrix}
5 & 0 \\
0 & 1
\end{matrix}
\right],
$$
which corresponds to principle curvatures $5$ and $1$ in the coordinate directions.
This leads to the formation of two cylinders in each side of $\partial_D \Omega$, as observed in \cite{AlBSm}. Our numerical experiments exhibit self-crossing before the two cylinders separate this time. Therefore, the lack of self-crossing alluded to in Example \ref{S:clamped} is not generic.

\subsection{Free Plate: Anisotropic Curvature}\label{ss:free}

\begin{figure}[htb]
  \begin{center}
	\includegraphics[scale=0.6]{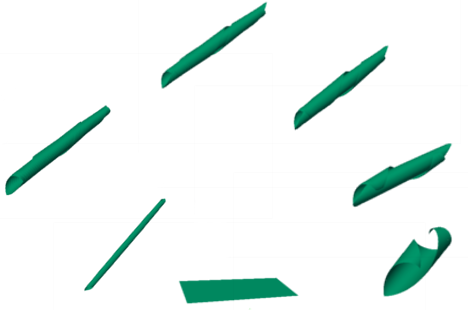}
  \end{center}
	\caption[Bilayer: Free Anisotropic]{Pseudo-evolution (counter-clockwise) towards equilibrium of a rectangular plate with anisotropic spontaneous curvature and free boundary conditions. The bilayer plate is depicted at times 0.0, 0.1, 0.5, 1.5, 2.5, 5.3, 22.0  $\times 10^{3}$ of the gradient flow. The last three snapshots reveal that the plate assumes a tighter configuration as an effort to become a full cylinder but instead rolls into a cigar.}
	\label{F:Anisotropic_Free}
\end{figure}

We now explore a cigar-type configuration motivated by work in mechanical engineering \cite{Simpson}. The plate \modif{$\Omega=(-5,5)\times(-2,2)$} is completely free of boundary conditions and has the same anisotropic spontaneous curvature
$$
Z =
\left[
\begin{matrix}
3 & -2 \\
-2 & 3
\end{matrix}
\right]
$$
as Example \ref{S:anisotropy}.
We observe in Figure\nobreakspace \ref {F:Anisotropic_Free} that the plate deforms at $45^o$ degrees with respect to the cartesian axes, similarly to Figure\nobreakspace \ref {F:Anisotropic} but, in the absence of a clamped side, it does so in a symmetric way and eventually reaches a cylindrical configuration (cigar).

\FloatBarrier
\FloatBarrier

\subsection{Free Plate: Wavy Pattern}

Several geometric shapes, including wave patterns, are obtained experimentally in \cite{JHZ} using bilayer and multilayer materials as building blocks for more complex self-folding and self-organizing structures. 
The wave pattern is created by alternating the position of the polymer that acts as the bilayer. 

\begin{figure}[htb]
  \begin{center}
	\includegraphics[scale=0.6]{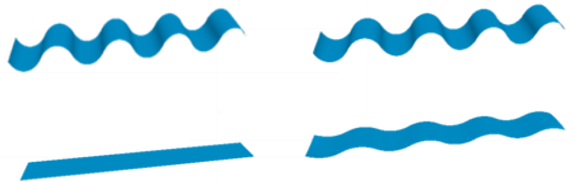}
  \end{center}
	\caption[Bilayer: Free Wave]{Pseudo-evolution (bottom-left to top-left) towards equilibrium of a rectangular plate with aspect ratio $8$ and alternating spontaneous curvature $Z=\pm I$ in each of its $8$ square parts. The bilayer plate is depicted at times 0.0, 0.1, 1.0, 2.3 $\times 10^{3}$ of the gradient flow.}
	\label{F:Wavy}
\end{figure}

We present here numerical simulations of wave patterns obtained by  
splitting the plate \modif{$\Omega=(-8,8)\times(-1,1)$} into eight equal parts along the $x_1$ direction. We alternate the spontaneous curvature between $Z=-I$ and $Z=I$
in each part, and impose free boundary conditions on the plate. 
The relaxation dynamics towards equilibrium is illustrated in Figure\nobreakspace \ref{F:Wavy}.
We note that the gradient flow is much faster than previous examples and, as a consequence, requires much fewer steps: the plate reaches the fourth depicted configuration quickly, which is then succeeded by  very small variations of shape and energy.

\FloatBarrier
\FloatBarrier

\FloatBarrier
\begin{figure}[htb]
  \begin{center}
    \includegraphics[scale=0.42]{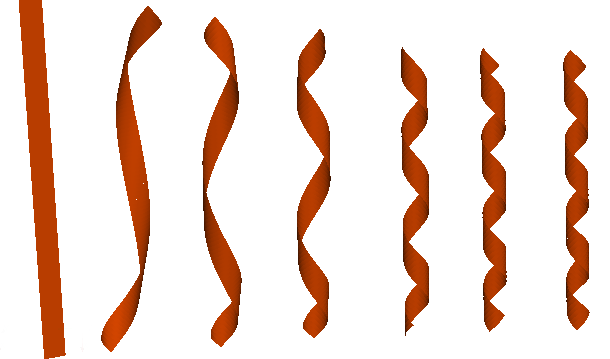}
  \end{center}
    \caption[Bilayer: DNA]{Pseudo-evolution (left to right) towards equilibrium of a rectangular plate with high aspect ratio and anisotropic spontaneous curvature with principal directions forming an angle of $\pi/4$ with the axes and principal curvatures $-1/2$ and $5/2$. The bilayer plate is depicted at times 0.0, 0.5, 2.0, 5.0, 15.0, 48.0 $\times 10^{3}$ of the gradient flow.}
	\label{F:DNA}
\end{figure} 
\FloatBarrier

\subsection{Free Plate: Helix Shape}\label{ss:DNA}

We now present the second example from \cite{JHZ}, which is a 
DNA-like shape.
We consider a high aspect ratio plate \modif{$\Omega= (-8,8) \times (-0.5,0.5)$,} which deforms with free boundary conditions under the effect of the spontaneous curvature 
$$
Z =
\left[
\begin{matrix}
1 & -3/2 \\
-3/2 & 1
\end{matrix}
\right].
$$
This choice of spontaneous curvature corresponds to principal directions that form an angle of $45$ degrees with the coordinate axes, similar to Example \ref{ss:free} (see Figure \ref{F:Anisotropic_Free}), but with principal curvatures $-1/2$ and $5/2$. We observe numerically that the relative magnitude of principal curvatures and aspect ratio leads to a deformation that resembles the twisting of DNA molecules (DNA-like helix). We illustrate the resulting deformation in \Cref{F:DNA} for several instances of the gradient flow.  Different combinations would lead to deformations of similar nature, but with different visual results in terms of plate twisting.

\section{Conclusions}\label{S:conclusions}

  In \cite{BaBoNo} we introduce a model for bilayer plates that undergo large
  (geometrically nonlinear) isometric deformations driven by an intrinsic spontaneous
  curvature tensor; see \cite{BasAbLaGr,JaSmIn,KuLPL,SchmidtEb,SmInLu} for enginnering applications and \cite{Schmidt} for analysis.
  In \cite{BaBoMuNo} we discuss a simplified model for thermal actuation of such
  bilayer plates. In both cases, we discretize the models with Kirchhoff elements
  and prove $\Gamma$-convergence in \cite{BaBoNo}.
  In \cite{BoNoNt} we propose a discontinuous Galerkin method (dG) for single layer plates that also display large isometric deformations \cite{Bartels,FJM}. The discrete energy functional captures the discontinuities of the discrete space and $\Gamma-$converges to the continuous energy \cite{BoNoNt}. In this paper we extend the dG methodology of \cite{BoNoNt} to the simplified bilayer model of \cite{BaBoMuNo}. Our contributions and pending questions are the following:
  \begin{enumerate}[$\bullet$]
    \item
      {\it Computational modeling}:
  The reduced model of \cite{BaBoMuNo} incorporates the effect of spontaneous curvature via a simple additional term to the elastic energy of \cite{Bartels,BoNoNt}. We propose a discrete gradient flow that treats such term explicitly and without interelement jumps. This, together with a linearized isometry constraint, leads to a linear saddle point problem to be solved at each step with a Schur complement algorithm. The latter exploits that the matrix for the inner loop is independent of the step and can thus be factored out only once; the corresponding linear system is solved with a parallel direct method. 

\smallskip
    \item
      {\it Simulations}:
      We showcase the excellent computational performance of dG with examples extracted from \cite{BaBoMuNo,BaBoNo} as well as from the engineering literature \cite{AlBSm,JHZ,Simpson}. The dG method is able to reproduce configurations of interest in engineering and medical applications. Moreover, dG exhibits a higher geometric flexibility and accuracy than the Kirchhoff element approach of \cite{Bartels,BaBoMuNo,BaBoNo}; the current simulations solidify further the merits of dG already discovered in \cite{BoNoNt}. The relaxation dynamics is a mathematical devise to reach equilibrium but does not prevent self-crossing; see Examples \ref{S:clamped} and \ref{S:middle} as well as earlier computations \cite{Bartels,BaBoMuNo,BaBoNo}. Avoiding self-intersection within a physically meaningful dynamics for plates remains open; this question is addressed in \cite{Bartels-HB} for rods. Simulations were carried out within the software platform {\it deal.ii} \cite{dealii85,dealiipaper}.

\smallskip      
    \item
      {\it Boundary conditions}:
We implement both Dirichlet and free boundary conditions with the Nitsche's approach, which enforces them weakly within the dG variational formulation rather than in the discrete space. This flexibility comes at the expense of a subtle dependence between the penalty parameters and the magnitudes of the initial energy and isometry defect. The choice of penalty parameters is thus a critical component of the dG approach. We are currently exploring parameter free options within the dG framework. 

\smallskip
    \item
      {\it Discrete gradient flow}:
 This relaxation dynamics is a semi-implicit discretization of a continuous gradient flow that hinges on orthogonality properties valid at the continuous level for isometries. Since such properties are violated slightly at the discrete level, it is no longer obvious that the discrete gradient flow decreases the discrete energy and guarantees the isometry defect. These properties were crucial in \cite{BoNoNt} to prove $\Gamma$-convergence.

\smallskip 
    \item
      {\it Future research}: There are two important questions that require further attention and remain open. The first is the numerical analysis of the current dG scheme. This involves the study of the discrete gradient flow \eqref{E:grad-flow2} and the $\Gamma$-convergence, including the analysis of the boundary-free case. As we mention after \eqref{E:energynormN}, coercivity for the boundary-free case holds only in the subspace of our discrete space of functions having vanishing zero and first moments. Elucidating how the discrete gradient flow affects this property is a crucial aspect of this analysis. The second question is the potential extension of our method to the physically interesting and challenging prestrained models \cite{Lewicka}.

  \end{enumerate}

\bibliographystyle{amsplain}

\end{document}